\documentclass{amsart}
\usepackage[all]{xy}
\usepackage{amsmath}
\usepackage{latexsym}
\RequirePackage{amssymb}
\RequirePackage[T1]{fontenc}

\newtheorem{dfn}{Definition}[section]
\newtheorem{prop}[dfn]{Proposition}

\newtheorem{qtn}[dfn]{Question}
\newtheorem{cjt}[dfn]{Conjecture}
\newtheorem{thm}[dfn]{Theorem}
\newtheorem{lem}[dfn]{Lemma}
\newtheorem{ex}[dfn]{Example}

\newtheorem{rmq}[dfn]{Remark}

\newcommand{\R}{\mathbb{R}}

\newcommand{\C}{\mathbb{C}}
\newcommand{\K}{\textbf{k}}

\newcommand{\pa}{\partial}
\newcommand{\al}{\mathfrak g}
\newcommand{\hl}{\mathfrak h}
\newcommand{\ml}{\mathfrak m}
\newcommand{\gun}{\mathfrak g_1}
\newcommand{\gde}{\mathfrak g_2}
\newcommand{\Gun}{\mathcal G_1}
\newcommand{\Gde}{\mathcal G_2}
\newcommand{\ggun}{C^c(\mathfrak g_1[2])}
\newcommand{\ggde}{C^c(\mathfrak g_2[2])}

\newcommand{\cyb}{\underline{{\rm CDYB}}}
\newcommand{\dt}{\underline{{\rm ADT}}}

\author{Damien Calaque}

\title[Quantization of formal classical dynamical $r$-matrices]
{Quantization of formal classical dynamical $r$-matrices: the reductive case}

\begin{document}

\maketitle

\begin{abstract}
In this paper we prove the existence of a formal dynamical twist quantization for any triangular and non-modified formal 
classical dynamical $r$-matrix in the reductive case. The dynamical twist is constructed as the 
image of the dynamical $r$-matrix by a $L_\infty$-quasi-isomorphism. This quasi-isomorphism also allows 
us to classify formal dynamical twist quantizations up to gauge equivalence. 
\end{abstract}

\section*{Introduction}

In \cite{F}, Felder introduced dynamical versions of both classical and quantum Yang-Baxter equations which has been 
generalized to the case of a nonabelian base in \cite{EV} for the classical part and in \cite{X} for the quantum part. 
Naturally this leads to quantization problems which have been formulated in terms of twist quantization \emph{à la} 
Drinfeld (\cite{D}) in \cite{X1,X,EE1,EE2}. 

Let us formulate this problem in the general context. Consider an inclusion $\hl\subset\al$ of Lie algebras equipped with 
an element $Z\in(\wedge^3\al)^\al$. A \emph{(modified) classical dynamical $r$-matrix} for $(\al,\hl,Z)$ is a regular 
(meaning $C^\infty$, meromorphic, formal, $\dots$ depending on the context) $\hl$-equivariant map 
$\rho:\hl^*\to\wedge^2\al$ which satisfies the \emph{(modified) classical dynamical Yang-Baxter equation} (CDYBE)
\begin{equation}\label{eq:cdybe}
{\rm CYB}(\rho)-{\rm Alt}(d\rho)=Z
\end{equation}
where ${\rm CYB}(\rho):=[\rho^{1,2},\rho^{1,3}]+[\rho^{1,2},\rho^{2,3}]+[\rho^{1,3},\rho^{2,3}]=\frac12[\rho,\rho]$ and 
$${\rm Alt}(d\rho):=\sum_i\big(h_i^1\frac{\pa\rho^{2,3}}{\pa\lambda^i}-
h_i^2\frac{\pa\rho^{1,3}}{\pa\lambda^i}+h_i^3\frac{\pa\rho^{1,2}}{\pa\lambda^i}\big)$$
Here $(h_i)$ and $(\lambda^i)$ are dual basis of $\hl$ and $\hl^*$. 

Let $\Phi=1+O(\hbar^2)\in(U\al^{\otimes3})^\al[[\hbar]]$ be an associator quantizing $Z$ (of which the existence was 
proved in \cite[proposition 3.10]{D2}). 
A \emph{dynamical twist quantization} of a (modified) classical dynamical $r$-matrix $\rho$ associated to $\Phi$ is a 
regular $\hl$-equivariant map $J=1+O(\hbar)\in{\rm Reg}(\hl^*,U\al^{\otimes2})[[\hbar]]$ such that 
${\rm Alt}\frac{J-1}{\hbar}=\rho~{\rm mod}~\hbar$ and which satisfies the \emph{(modified) dynamical twist equation} (DTE)
\begin{equation}\label{eq:dte}
J^{12,3}(\lambda)*J^{1,2}(\lambda+\hbar h^3)=\Phi^{-1}J^{1,23}(\lambda)*J^{2,3}(\lambda)
\end{equation}
where $*$ denotes the PBW star-product of functions on $\hl^*$ and 
$$J^{1,2}(\lambda+\hbar h^3):=\sum_{k\geq0}\frac{\hbar^k}{k!}\sum_{i_1,\dots,i_k}
(\pa_{\lambda^{i_1}}\cdots\pa_{\lambda^{i_k}}J)(\lambda)\otimes(h_{i_1}\cdots h_{i_k})$$

Now observe that many (modified) classical dynamical $r$-matrices can be viewed as formal ones by taking their Taylor 
expansion at $0$. In this paper we are interested in the following conjecture: 
\begin{cjt}[\cite{EE1}]\label{thm:cjt}
Any (modified) formal classical dynamical $r$-matrix admits a dynamical twist quantization. 
\end{cjt}
Let us reformulate DTE in the formal framework. A formal (modified) dynamical twist is an element 
$J(\lambda)=1+O(\hbar)\in(U\al^{\otimes2}\hat\otimes\hat S\hl)^\hl[[\hbar]]$ which satisfies DTE, and 
$J^{1,2}(\lambda+\hbar h^3)\in(U\al^{\otimes3}\hat\otimes\hat S\hl)[[\hbar]]$ is equal to 
$({\rm id}^{\otimes2}\otimes\tilde\Delta)(J)$ where $\tilde\Delta:\hat S\hl\to(U\al\hat\otimes\hat S\hl)[[\hbar]]$ is 
induced by $\hl\ni x\mapsto\hbar x\otimes1+1\otimes x$. Then define 
$K:=J(\hbar\lambda)\in(U\al^{\otimes2}\otimes S\hl)^\hl[[\hbar]]$ 
which we view as an element of $(U\al^{\otimes2}\otimes U\hl)^\hl[[\hbar]]$ using the symmetrization map $S\hl\to U\hl$. 
Since $J$ is a solution of DTE $K$ satisfies the \emph{(modified) algebraic dynamical twist equation} (ADTE) 
\begin{equation}\label{eq:adte}
K^{12,3,4}K^{1,2,34}=(\Phi^{-1})^{1,2,3}K^{1,23,4}K^{2,3,4}
\end{equation}
Moreover and by construction, $K=1+\sum_{n\geq1}\hbar^nK_n$ has the \emph{$\hbar$-adic valuation property}. Namely, $U\hl$ 
is filtered by $(U\hl)_{\leq n}=\ker{({\rm id}-\eta\circ\varepsilon)^{\otimes n+1}\circ\Delta^{(n)}}$ where 
$\varepsilon:U\hl\to\K$ and $\eta:\K\to U\hl$ are the counit and unit maps, and $K_n\in(U\hl)_{\leq n-1}$. 
Conversely, any algebraic dynamical twist having the $\hbar$-adic valuation property can be obtained from a unique 
formal dynamical twist by this procedure. \\

This paper, in which we always assume $Z=0$ and $\Phi=1$ (non-modified case), is organized as follow. 

In section 1 we define two differential graded Lie algebras ({\tt dgla}'s) respectively associated to classical 
dynamical $r$-matrices and algebraic dynamical twists. Then we formulate the main theorem of this paper which states that 
if $\hl$ admits an ${\rm ad}\hl$-invariant complement (the \emph{reductive} case) then these two {\tt dgla}'s are 
$L_\infty$-quasi-isomorphic and we prove that it implies Conjecture \ref{thm:cjt} in this case, which generalizes 
Theorem 5.3 of \cite{X1}: 
\begin{thm}\label{thm:rmat}
In the reductive case, any formal classical dynamical $r$-matrix for $(\al,\hl,0)$ admits a dynamical twist 
quantization (associated to the trivial associator). 
\end{thm}
The second section is devoted to the proof of the main theorem of section 1: using an equivariant formality theorem for 
homogeneous spaces which is obtain from \cite{Do}, we construct a $L_\infty$-quasi-isomorphism which we then modify in 
order to obtain the desired one. 
We use this $L_\infty$-quasi-isomorphism to classify formal dynamical twist quantizations up to gauge equivalence for the 
reductive case in section 3. 
In section 4 we prove that if $\al=\hl\oplus\ml$ for $\hl$ abelian and $\ml$ a Lie subalgebra then the results of 
sections 1 and 2 are still true in this situation. 
We conclude the paper with some open questions, and recall basic results for $L_\infty$-algebras in an appendix. \\

\noindent{\bf Acknowledgements. }I thank Benjamin Enriquez for many usefull discussions on this subject. 

\section{Definitions and results}

Let $\hl\subset\al$ be an inclusion of Lie algebras. 

\subsection{Algebraic structures associated to CDYBE}\label{cyb}

Let us consider the following graded vector space 
$$\cyb:=\wedge^*\al\otimes S\hl=\bigoplus_{k\geq0}\wedge^{k}\al\otimes S\hl$$
equipped with the differential $\textrm d$ defined by 
\begin{equation}\label{eq:diff}
\textrm d(x_1\wedge\dots\wedge x_{k}\otimes h_1\cdots h_l)
:=-\sum_{i=1}^lh_i\wedge x_1\wedge\dots\wedge x_k\otimes h_1\cdots h_l\hat h_i
\end{equation}
With the exterior product $\wedge$ it becomes a differential graded commutative associative algebra. 
Moreover, one can define a graded Lie bracket of degree $-1$ on $\cyb$ which is the Lie bracket of $\al$ extended to 
$\cyb$ in the following way: 
\begin{equation}\label{eq:GA}
[a,b\wedge c]=[a,b]\wedge c+(-1)^{(| a|-1)| b|}b\wedge[a,c]
\end{equation}
Thus one can observe that polynomial solutions to CDYBE are exactly elements $\rho\in\cyb$ of degree $2$ such that 
$\textrm d\rho+\frac12[\rho,\rho]=0$. We would like to say that such a $\rho$ is a Maurer-Cartan element but 
$(\cyb[1],\textrm d,[,])$ is not a differential graded Lie algebra ({\tt dgla}). 

Instead, remember that we are interested in $\hl$-equivariant solutions of CDYBE (i.e., dynamical $r$-matrices) and 
thus consider the subspace $\gun=(\cyb)^{\hl}$ of $\hl$-invariants with the same differential and bracket. 
\begin{prop}
$(\gun[1],{\rm d},[,])$ is a {\tt dgla}. Moreover $(\gun,{\rm d},\wedge,[,])$ is a Gerstenhaber algebra. 
\end{prop}
\begin{proof}
Let $a=x_1\wedge\dots\wedge x_k\otimes h_1\cdots h_s$ and $b=y_1\wedge\dots\wedge y_l\otimes m_1\cdots m_t$ be 
$\hl$-invariant elements in $\gun$. We want to show that 
\begin{equation}\label{eq:dgla}
\textrm d[a,b]=[\textrm da,b]+(-1)^{k-1}[a,\textrm db]
\end{equation}
The l.h.s. of (\ref{eq:dgla}) is equal to 
\begin{eqnarray*}
-\Big(\sum_{i=1}^sh_i\wedge[x_1\wedge\dots\wedge x_k,y_1\wedge\dots\wedge y_l]\otimes h_1\cdots h_sm_1\cdots m_t\hat h_i
&&\\
+\sum_{j=1}^tm_j\wedge[x_1\wedge\dots\wedge x_k,y_1\wedge\dots\wedge y_l]\otimes h_1\cdots h_sm_1\cdots m_t\hat m_j\Big)&&
\end{eqnarray*}
The first term in the r.h.s. of (\ref{eq:dgla}) gives 
$$\sum_{i=1}^s\big((-1)^{k-1}x_1\wedge\cdots\wedge x_k\wedge[h_i,y_1\wedge\cdots\wedge y_l]
-h_i\wedge[x_1\wedge\dots\wedge x_k,y_1\wedge\dots\wedge y_l]\big)\otimes h_1\cdots h_sm_1\cdots m_t\hat h_i$$
and for the second term we obtain 
$$\sum_{j=1}^t\big((-1)^{k-1}[m_j,x_1\wedge\cdots\wedge x_k]\wedge y_1\wedge\cdots\wedge y_l
-m_j\wedge[x_1\wedge\dots\wedge x_k,y_1\wedge\dots\wedge y_l]\big)\otimes h_1\cdots h_sm_1\cdots m_t\hat m_j$$
Thus the difference between the l.h.s. and the r.h.s. of (\ref{eq:dgla}) is equal to 
\begin{eqnarray*}
(-1)^k\Big(
\sum_{i=1}^kx_1\wedge\cdots\wedge x_k\wedge[h_i,y_1\wedge\cdots\wedge y_l]\otimes h_1\cdots h_sm_1\cdots m_t\hat h_i &&\\
+\sum_{j=1}^l[m_j,x_1\wedge\cdots\wedge x_k]\wedge y_1\wedge\cdots\wedge y_l\otimes h_1\cdots h_sm_1\cdots m_t\hat m_j\Big) 
&&
\end{eqnarray*}
Then using $\hl$-invariance of $a$ and $b$ one obtains 
$$(-1)^{k-1}\sum_{i,j}x_1\wedge\cdots\wedge x_k\wedge y_1\wedge\cdots\wedge y_l\otimes
\big(h_1\cdots h_sm_1\cdots m_t([h_i,m_j]-[m_j,h_i])\hat h_i\hat m_j\big)=0$$

The second statement of the proposition is obvious from the definition (\ref{eq:GA}) of the bracket. 
\end{proof}
Let $\rho(\lambda)\in(\wedge^2\al\hat\otimes\hat S\hl)^\hl$ be a formal classical dynamical $r$-matrix. Since $\rho$ 
satisfies CDYBE, $\alpha:=\hbar\rho(\hbar\lambda)\in\hbar\gun[[\hbar]]$ is a Maurer-Cartan element ({\rm i.e.} 
${\rm d}\alpha+\frac12[\alpha,\alpha]=0$). 

\subsection{Algebraic structures associated to ADTE}

Let us now consider the graded vector space 
$$\dt:=T^*U\al\otimes U\hl=\bigoplus_{k\geq0}\otimes^kU\al\otimes U\hl$$
equipped with the differential $b$ given by 
\begin{equation}\label{eq:b}
b(P):=P^{2,\dots,k+2}+\sum_{i=1}^{k+1}(-1)^iP^{1,\dots,ii+1,\dots,k+2}
\quad{\rm for}\quad P\in\otimes^kU\al\otimes U\hl
\end{equation}
\begin{rmq}\emph{This is just the coboundary operator of Hochschild's cohomology with value in a comodule; and $b^2=0$ follows 
directly from an easy calculation. }
\end{rmq}

One can define on $\dt$ an associative product $\cup$ (the \emph{cup} product) which is given on homogeneous elements 
$P\in\otimes^{k}U\al\otimes U\hl$ and $Q\in\otimes^{l}U\al\otimes U\hl$ by 
$$P\cup Q:=P^{1,\dots,k,k+1\dots k+l+1}Q^{k+1,\dots,k+l+1}$$
\begin{prop}
$(\dt,b,\cup)$ is a differential graded associative algebra. 
\end{prop}
\begin{proof}
The cup product is obviously associative. Thus the only thing we have to check is that 
\begin{equation}\label{eq:dgaa}
b(P\cup Q)=bP\cup Q+(-1)^{| P|}P\cup bQ
\end{equation}
Let $k=| P|$ and $l=| Q|$. The l.h.s. of (\ref{eq:dgaa}) is equal to 
\begin{eqnarray*}
P^{2,\dots,k+1,k+2\dots k+l+2}Q^{k+2,\dots,k+l+2}
+\sum_{i=1}^{k}(-1)^iP^{1,\dots,ii+1,\dots,k+1,k+2\dots k+l+2}Q^{k+2,\dots,k+l+2} && \\
+\sum_{i=k+1}^{k+l+1}(-1)^iP^{1,\dots,k,k+1\dots k+l+2}Q^{k+1,\dots,ii+1,\dots,k+l+2}
\end{eqnarray*}
The first line of this expression is equal to 
$$bP\cup Q-(-1)^{k+1}P^{1,\dots,k,k+1\dots k+l+2}Q^{k+2,\dots,k+l+2}$$
and the last term of the same expression gives 
$$(-1)^k\big(P\cup bQ-P^{1,\dots,k,k+1\dots k+l+2}Q^{k+2,\dots,k+l+2}\big)$$
The proposition is proved. 
\end{proof}

Recall that in the case $\hl=\{0\}$ one can define a brace algebra structure on $(T^*U\al)[1]$ (see \cite{G}). 
Unfortunately we are not able to extend this structure to $\dt$ in general. Since we deal with $\hl$-equivariant solutions 
of ADTE we can consider the subspace $\gde=(\dt)^\hl$ of $\hl$-invariants. Let us now define a collection of linear 
homogeneous maps of degree zero $\{-|-,\dots,-\}:\gde[1]\otimes\gde[1]^{\otimes m}\to\gde[1]$ indexed by $m\geq0$, and 
$\{P|Q_1,\dots,Q_m\}$ is given by 
$$\sum_{\substack{0\leq i_1,i_m+k_m\leq n\\i_l+k_l\leq i_{l+1}}}
(-1)^\epsilon P^{1,\dots,i_1+1\dots i_1+k_1,\dots,i_m+1\dots i_m+k_m,\dots,n+1}
\prod_{s=i}^mQ_s^{i_s+1,\dots,i_s+k_s,i_s+k_s+1\dots n+1}$$
where $k_s=|Q_s|$, $n=|P|+\sum_sk_s-m$ and $\epsilon=\sum_s(k_s-1)i_s$. 
\begin{prop}
$(\gde[1],\{-|-,\dots,-\})$ is a brace algebra. 
\end{prop}
\begin{proof}
Since we work with $\hl$-invariant elements one can remark that if $i_s+k_s\leq i_t$ then 
$Q_s^{i_s+1,\dots,i_s+k_s,i_s+k_s+1\dots n+1}$ and $Q_t^{i_t+1,\dots,i_t+k_t,i_t+k_t+1\dots n+1}$ commute. 
Using this the proof becomes identical to the case when $\hl=0$ (see \cite{G} for example). 
\end{proof}
Now observe that since $m=1^{\otimes3}\in(\otimes^2U\al\otimes U\hl)^\hl$ is such that $\{m|m\}=0$ one obtains a 
$B_\infty$-algebra structure (\cite{B}) on $\gde$ (see \cite{Kh}). More precisely, we have a differential graded bialgebra 
structure on the cofree tensorial coalgebra $T(\gde[1])$ of which structure maps $a^n,a^{p,q}$ are given by
\begin{itemize}
\item $a^1(P)=bP=(-1)^{|P|-1}[m,P]_G$, where 
$$[P,Q]_G:=\{P|Q\}-(-1)^{(|P|-1)(|Q|-1)}\{Q|P\}$$
\item $a^2(P,Q)=\{m|P,Q\}=P\cup Q$
\item $a^{0,1}=a^{1,0}={\rm id}$
\item $a^{1,n}(P;Q_1,\dots,Q_n)=\{P|Q_1,\dots,Q_n\}$ for $n\geq1$
\item all other maps are zero
\end{itemize}
In particular, we have 
\begin{prop}
$(\gde[1],b,[,]_G)$ is a {\tt dgla}. 
\end{prop}
\begin{rmq}\emph{Since that for any graded vector space $V$, {\tt dg} bialgebra structures on the cofree 
coassociative coalgebra $T^cV$ are in one-to-one correspondence with {\tt dg} Lie bialgebra structures on the cofree Lie 
coalgebra $L^cV$ (see \cite{T}, section 5), then $L^c(\gde[1])$ becomes a {\tt dg} Lie bialgebra with differential and Lie 
bracket given by maps $l^n,l^{p,q}$ such that $l^1=b$ and $l^{1,1}=[,]_G$. Therefore 
$d_2:=\sum_{i\geq0}l^i+\sum_{p,q\geq0}l^{p,q}:C^c(L^c(\mathfrak g_2[1]))\to C^c(L^c(\mathfrak g_2[1]))$ 
defines a $G_\infty$-algebra structure on $\gde$ ($d_2\circ d_2=0$ since $d_2$ is just the Chevalley-Eilenberg 
differential on the {\tt dg} Lie algebra $L^c(\gde[1])$). }
\end{rmq}

\subsection{Main result and proof of theorem \ref{thm:rmat}}\label{sec:proof}

First of all, observe that $\cyb$, $\gun$ and $\Gun:=\ggun$ have a natural grading induced by the one of $S\hl$. In the 
same way $\dt$, $\gde$ and $\Gde:=\ggde$ have a natural filtration induced by the one of $U\hl$. 
Our main goal is to prove the following theorem, which is sufficient to obtain algebraic dynamical twists from formal 
dynamical $r$-matrices. 
\begin{thm}\label{thm:main}
In the reductive case, there exists a $L_\infty$-quasi-isomorphism 
$$\Psi:(\Gun,{\rm d}+[,])\to(\Gde,b+[,]_G)$$
with the following two filtration properties: 
\begin{enumerate}
\item[(F1)] $\forall X\in(\gun)_k$, $\Psi^1(X)=({\rm alt}\otimes{\rm sym})(X)~{\rm mod}~(\gde)_{\leq k-1}$
\item[(F2)] $\forall X\in(\Lambda^n\gun)_k$, $\Psi^n(X)\in(\gde)_{\leq n+k-1}$
\end{enumerate}
\end{thm}
\begin{proof}[Proof of Theorem \ref{thm:rmat}]
Now consider a formal solution $\rho(\lambda)\in(\wedge^2\al\hat\otimes\hat S\hl)^\hl$ to CDYBE. Let us define 
$\alpha:=\hbar\rho(\hbar\lambda)\in\hbar\gun[[\hbar]]$ which is a Maurer-Cartan element in $\hbar\gun[[\hbar]]$. The 
$L_\infty$-morphism property implies that $\widetilde\alpha:=\sum_{n=1}^\infty\frac1{n!}\Psi^n(\Lambda^n\alpha)$ is a 
Maurer-Cartan element in $\hbar\gde[[\hbar]]$; this exactly means that 
$K:=1+\widetilde\alpha\in(\otimes^2U\al\otimes U\hl)^\hl[[\hbar]]$ satisfies ADTE. Moreover, due to (F2) the coefficient 
$K_n$ of $\hbar^n$ in $K$ lies in $(\gde)_{\leq n-1}$. It means that there exists 
$J\in(U\al^{\otimes 2}\hat\otimes\hat S\hl)^\hl[[\hbar]]$ satisfying DTE and such that 
$K=({\rm id}^{\otimes2}\otimes{\rm sym})(J(\hbar\lambda))$. Finally, property (F1) obviously implies that the semi-classical 
limit condition $\frac{J-J^{{\rm op}}}{\hbar}=\rho\textrm{ mod }\hbar$ is satisfied. 
\end{proof}

\section{Proof of theorem \ref{thm:main}}

In this section we assume that $\al=\hl\oplus\ml$ with $[\hl,\ml]\subset\ml$. Let us denote by ${\rm p}:\al\rightarrow\ml$ 
the projection on $\ml$ along $\hl$; it is $\hl$-equivariant. 

\subsection{Resolutions}\label{sec:reso}

Let us first observe that the bilinear map $[,]_\ml:=(\wedge^{\cdot}{\rm p})\circ[,]$ defines a graded Lie bracket of 
degree $-1$ on $(\wedge^*\ml)^\hl$. Then we prove 
\begin{prop}\label{reso1}
The natural map $p_1:(\gun[1],{\rm d},[,])\to((\wedge^*\ml)^\hl[1],0,[,]_\ml)$ is a morphism of {\tt dgla}'s. 
Moreover, there exists an operator $\delta:\gun^*\to\gun^{*-1}$ such that $\delta{\rm d}+{\rm d}\delta={\rm id}-p_1$, 
$\delta\circ\delta=0$ and $\delta\big((\gun)_k\big)\subset(\gun)_{k+1}$. In particular, $p_1$ induces an isomorphism in 
cohomology. 
\end{prop}
\begin{proof}
The projection $p_1:=(\wedge^{\cdot}{\rm p})\otimes\varepsilon:(\cyb,{\rm d})\to(\wedge^*\ml,0)$ is a $\hl$-equivariant 
morphism of complexes, and it obviously restricts to a morphism of (differential) graded Lie algebras at the level of 
$\hl$-invariants. 

Moreover, $\wedge^n\al\otimes S\hl\cong\bigoplus_{p+q=n}\wedge^p\ml\otimes\wedge^q\hl\otimes S\hl$ as a $\hl$-module; and 
under this identification ${\rm d}$ becomes $-{\rm id}\otimes d_K$, where 
$d_K:\wedge^*\hl\otimes S\hl\to\wedge^{*+1}\hl\otimes S\hl$ is Koszul's coboundary operator, and $p_1$ corresponds to the 
projection on the part of zero antisymmetric and symmetric degrees in $\hl$. Let us define $\delta={\rm id}\otimes\delta_K$ 
with $\delta_K:\wedge^*\hl\otimes S^*\hl\to\wedge^{*-1}\hl\otimes S^{*+1}\hl$ defined by 
$$\delta_K(x_1\wedge\cdots\wedge x_n\otimes h_1\cdots h_m)=
\left\{{\small
\begin{array}{cl}
\frac{1}{m+n}\sum_i(-1)^ix_1\wedge\cdots\hat{x_i}\cdots\wedge x_n\otimes h_1\cdots h_mx_i
& \mbox{if $m+n\neq0$} \\
0 & \mbox{otherwise}
\end{array}}
\right.$$
Finally remark $\delta$ is a $\hl$-equivariant homotopy operator: $\delta{\rm d}+{\rm d}\delta={\rm id}-p_1$ and 
$\delta\circ\delta=0$. The proposition is proved. 
\end{proof}
Now we prove a similar result for $\gde$. 
Let us first define $U\ml:={\rm sym}(S\ml)\subset U\al$; this is a sub-coalgebra of $U\al$ and thus $T^*U\ml$ equipped 
with its Hochschild's coboundary operator $b_\ml$ becomes a cochain subcomplex of the Hochschild complex 
$(T^*U\al,b_\al)$ of $U\al$. We also have the following 
\begin{lem}\label{lem1}
$U\al=U\al\cdot\hl\oplus U\ml$ as a filtered $\hl$-module. Moreover $[,]_{G,\ml}:=(\otimes^{\cdot}{\rm p})\circ[,]$ 
defines a graded Lie bracket of degree $-1$ on $(T^*U\ml)^\hl$
\end{lem}
\begin{proof}
See \cite[Ch.II \S4.2]{H} for the first statement. The second statement follows from a direct computation. 
\end{proof}
Then we prove the 
\begin{prop}\label{reso2}
The natural map $p_2:(\gde[1],b,[,]_G)\to((T^*U\ml)^\hl[1],b_\ml,[,]_{G,\ml})$ is a morphism of {\tt dgla}'s. 
Moreover, there exists an operator $\kappa:\gde^*\to\gde^{*-1}$ such that $\kappa b+b\kappa={\rm id}-p_2$, 
$\kappa\circ\kappa=0$ and $\kappa\big((\gde)_{\leq k}\big)\subset(\gde)_{\leq k+1}$. In particular, $p_2$ induces an 
isomorphism in cohomology. 
\end{prop}
\begin{proof}
The projection $p_2:=(\otimes^{\cdot}{\rm p})\otimes\varepsilon:(\dt,b)\to(T^*U\ml,b_\ml)$ is a $\hl$-equivariant 
morphism of complexes, and it obviously restricts to a morphism of {\tt dgla}'s at the level of $\hl$-invariants (by 
lemma \ref{lem1}). 

Remember that $\gde$ has a natural filtration induced by the one of $U\hl$. Then one obtains a spectral sequence of which 
we compute the first terms: 
$$\begin{array}{lc}
E_0^{*,*}=(T^*U\al\otimes S^*\hl)^\hl & d_0=b_\al\otimes{\rm id} \\
E_1^{*,*}=(\wedge^*\al\otimes S^*\hl)^\hl & d_1={\rm d} \\
E_2^{*,*}=E_2^{*,0}=(\wedge^*\ml)^\hl & d_2=0
\end{array}$$
Then the proposition follows from proposition \ref{reso1}. 
\end{proof}

\subsection{Inverting $p_2$}\label{sec:p2}

In this subsection, taking our inspiration from \cite[appendix]{Mo}, we prove the following 
\begin{prop}
There exists a $L_\infty$-quasi-isomorphism 
$$\mathcal Q_2:(C^c((T^*U\ml)^\hl[2]),b_\ml+[,]_{G,\ml})\to(\ggde,b+[,]_G)$$
such that $\mathcal Q_2^1$ is the natural inclusion and $\mathcal Q_2^n$ takes values in $(\gde)_{\leq n-1}$. 
\end{prop}
\begin{proof}
Let $(N,b_N)\subset(\gde,b)$ be the kernel of the surjective morphism of complexes $p_2:(\gde,b)\to((T^*U\ml)^\hl,b_\ml)$. 
It follows from the proofs of propositions \ref{reso1} and \ref{reso2} that there exists an operator $H:N^*\to N^{*-1}$ 
such that $H\circ H=0$, $b_NH+Hb_N={\rm id}$ and $H(N_{\leq n})\subset N_{\leq n+1}$. 

Now let us construct a $L_\infty$-isomorphism 
$$\mathcal F:\big(\ggde,b+[,]_G\big)\tilde\longrightarrow\big(C^c((T^*U\ml)^\hl[2]\oplus N[2]),b_\ml+b_N+[,]_{G,\ml}\big)$$
with structure maps $\mathcal F^n:\Lambda^n\gde\to((T^*U\ml)^\hl\oplus N)[1-n]$ such that 
\begin{itemize}
\item $\mathcal F^1$ is the sum of $p_2$ with the projection on $N$ along $(T^*U\ml)^\hl$ (in some sense $\mathcal F^1$ is 
the identity), 
\item for any $n>1$ and $X\in(\Lambda^n\gde)_{\leq k}$, $\mathcal F^n(X)\in N_{\leq n+k-1}$. 
\end{itemize}
Let us prove it by induction on $n$. First $\mathcal F^1$ is a morphism of complexes by definition. Then let us define 
$\mathcal K_2:\Lambda^2\gde\to((T^*U\ml)^\hl\oplus N)[1]$ by 
$$\mathcal K_2(x\Lambda y)=[\mathcal F^1(x),\mathcal F^1(y)]_{G,\ml}-\mathcal F^1([x,y]_G)$$
It takes values in $N[1]$ and is such that $b_N\mathcal K_2(x,y)+\mathcal K_2(bx,y)+\mathcal K_2(x,by)=0$. Consequently 
$\mathcal F^2:=H\circ\mathcal K_2:\Lambda^2\gde\to N$ is such that 
$$b_N\mathcal F^2(x,y)-\mathcal F^2(bx,y)-\mathcal F^2(x,by)=\mathcal K_2(x,y)
\quad(L_\infty\textrm{-condition for }\mathcal F^2)$$
and for any $X\in(\Lambda^2\gde)_{\leq k}$, $\mathcal F^2(X)\in N_{\leq k+1}$. 
After this, suppose we have constructed $\mathcal F^1,\dots,\mathcal F^n$ and let us define 
$$\mathcal K_{n+1}:=[,]_{G,\ml}\circ\mathcal F^{\leq n}-\mathcal F^{\leq n}\circ[,]_G
:\Lambda^2\gde\to((T^*U\ml)^\hl\oplus N)[1]$$
It obviously takes values in $N[1]$ and is such that $b_NK_{n+1}+K_{n+1}b=0$. 
Consequently $\mathcal F^{n+1}:=H\circ K_{n+1}$ satisfies the $L_\infty$-condition 
$$b_N\mathcal F^{n+1}-\mathcal F^{n+1}b=b_NHK_{n+1}-HK_{n+1}b=(b_NH+Hb_N)K_{n+1}=K_{n+1}$$
and for any $X\in(\Lambda^n\gde)_{\leq n+1}$, $\mathcal F^{n+1}(X)\in N_{\leq n+k}$ 
(since $\mathcal K_{n+1}(X)\in N_{\leq n+k-1}$). 

Now let $\mathcal H$ be the inverse of the isomorphism $\mathcal F$, it is such that for any $n\geq1$ and 
$X\in(\Lambda^n\gde)_{\leq k}$, $\mathcal H^n(X)\in N_{\leq n+k-1}$. Finally we obtain $\mathcal Q_2$ by composing 
$\mathcal H$ with the inclusion of {\tt dgla}'s $(T^*U\ml)^\hl[1]\hookrightarrow((T^*U\ml)^\hl\oplus N)[1]$. 
\end{proof}

\subsection{End of the proof}\label{sec:end}

Recall from \cite[Ch.II \S4.2]{H} that $(T^*U\ml)^\hl={\rm Diff}^*(G/H)^G$ and 
$(\wedge^*\ml)^\hl=\Gamma(G/H,\wedge^*T(G/H))^G$ as {\tt dgla}'s. 
Remember also from \cite[Ch.II \S8]{N} that $G$-invariant connections on $G/H$ are in one-to-one correspondence with 
$\hl$-equivariant linear maps $\alpha:\ml\otimes\ml\to\ml$, and that the torsion tensor is given by 
$\alpha-\alpha^{21}-{\rm p}\circ[,]$. Thus $G/H$ is equipped with a $G$-invariant torsion free connection 
$\nabla$, corresponding to the map $\alpha:=\frac12{\rm p}\circ[,]$. 
Then using a theorem of Dolgushev, see \cite[theorem 5]{Do}, we obtain a $G$-equivariant $L_\infty$-quasi-isomorphism 
$\phi:\Gamma(G/H,\wedge^*T(G/H))\to{\rm Diff}^*(G/H)$ with first structure map $\phi^1={\rm alt}$, which restricts to a 
$L_\infty$-quasi-isomorphism at the level of $G$-invariants. 
Let us define $\psi:=\mathcal Q_2\circ\phi\circ p_1:(\ggun,{\rm d}+[,])\to(\ggde,b+[,]_G)$; it is a 
$L_\infty$-quasi-isomorphism with first structure map 
$\psi^1=({\rm alt}\otimes1)\circ(\wedge^{\cdot}{\rm p}\otimes\varepsilon)$. 

Finally define $V:=({\rm alt}\otimes{\rm sym})\circ\delta:\gun\to\gde[-1]$ and use lemma \ref{lem} to construct a 
$L_\infty$-quasi-morphism $\Psi:(\ggun,{\rm d}+[,])\to(\ggde,b+[,]_G)$ with first structure map 
$\Psi^1=\psi^1+b\circ V+V\circ{\rm d}$. Since for any $X\in(\Gun)_k$, then
$$\begin{array}{c}
b\circ({\rm alt}\otimes{\rm sym})(X)=({\rm alt}\otimes{\rm sym})\circ{\rm d}(X)\textrm{ mod }(\gde)_{\leq k-1} \\ \\
\begin{array}{rcl}
\textrm{ and thus }\quad\Psi^1(X) & = & \psi^1(X)+bV(X)+V({\rm d}X) \\
& = & ({\rm alt}\otimes{\rm sym})\circ(p_1+{\rm d}\delta+\delta{\rm d})(X)~{\rm mod}~(\gde)_{\leq k-1} \\
& = & ({\rm alt}\otimes{\rm sym})(X)~{\rm mod}~(\gde)_{\leq k-1}
\end{array}\end{array}$$
Consequently $\Psi$ satisfies (F1). Moreover, it follows from remark \ref{remark} that $\Psi$ also satisfies 
(F2). $\Box$

\section{Classification}

Theorem \ref{thm:main} implies a stronger result than just the existence of the twist quantization. Namely, since $\Psi$ 
is a $L_\infty$-quasi-isomorphism there is a bijection between the moduli spaces of Maurer-Cartan elements of {\tt dgla}'s 
$(\gun[1])[[\hbar]]$ and $(\gde[1])[[\hbar]]$.  

\subsection{Classification of algebraic and formal dynamical twists}

Following \cite{EE1}, two dynamical twists $J(\lambda)$ and $J'(\lambda)$ are said to be 
\emph{gauge equivalent} if there exists a regular $\hl$-equivariant map 
$T(\lambda)=\textrm{exp}(q)+O(\hbar)\in\textrm{Reg}(\hl^*,U\al)^\hl[[\hbar]]$, with 
$q\in\textrm{Reg}(\hl^*,\al)^\hl$ such that $q(0)=0$, and satisfying 
\begin{equation}\label{eq:gaugequiv}
J'(\lambda)=T^{12}(\lambda)*J(\lambda)*T^2(\lambda)^{-1}*T^1(\lambda+\hbar h^2)^{-1}
\end{equation}
Dealing with formal functions one can easily derive an equivalence relation for the corresponding 
algebraic dynamical twists $K=J(\hbar\lambda)$ and $K'=J'(\hbar\lambda)$: 
\begin{equation}\label{eq:algequiv}
K'=Q^{12,3}K(Q^{2,3})^{-1}(Q^{1,23})^{-1}
\end{equation}
in $(U\al^{\otimes2}\otimes U\hl)^\hl[[\hbar]]$, with $Q=1+O(\hbar)\in(U\al\otimes U\hl)^\hl[[\hbar]]$ given by 
$T(\hbar\lambda)$. \\

Assume now we are in the reductive case. 

Since the composition $\mathcal Q_2\circ\phi:(C^c((\wedge\ml)^\hl[2]),[,]_\ml)\to(C^c(\gde[2]),b+[,]_G)$ in the previous 
section is a $L_\infty$-quasi-isomorphism then we have a bijective correspondance 
\begin{equation}\label{eq:classiftwist}
\frac{\{\pi\in\hbar(\wedge^2\ml)^\hl[[\hbar]]~\textrm{s.t.}~[\pi,\pi]_\ml=0\}}{G_0}
\longleftrightarrow
\frac{\{\textrm{algebraic dynamical twists}\}}{\textrm{gauge equivalence }(\ref{eq:algequiv})}
\end{equation}
where $G_0$ is the prounipotent group corresponding to the Lie algebra $\hbar\,\ml^\hl[[\hbar]]$. 
Moreover, since the structure maps $\mathcal Q_2^n$ take values in $(\gde)_{\leq n-1}$ then it appears that any algebraic 
dynamical twist is gauge equivalent to a one with the $\hbar$-adic valuation property and thus we have a bijection 
\begin{equation}\label{eq:algform}
\frac{\{\textrm{algebraic dynamical twists}\}}{\textrm{gauge equivalence }(\ref{eq:algequiv})}
\longleftrightarrow
\frac{\{\textrm{formal dynamical twists}\}}{\textrm{gauge equivalence }(\ref{eq:gaugequiv})}
\end{equation}

\subsection{Classical counterpart}

Assume that we are in the reductive case. Since $p_1$ is a $L_\infty$-quasi-isomorphism by propostion 
\ref{reso1} then we have a bijection 
\begin{equation*}
\frac{\{\alpha\in\hbar(\wedge^2\al\otimes S\hl)^\hl[[\hbar]]~\textrm{s.t.}~{\rm d}\alpha+\frac12[\alpha,\alpha]=0\}}{G_1}
\longleftrightarrow
\frac{\{\pi\in\hbar(\wedge^2\ml)^\hl[[\hbar]]~\textrm{s.t.}~[\pi,\pi]_\ml=0\}}{G_0}
\end{equation*}
where $G_1$ is a prounipotent group and its action (by affine transformations) is given by the exponentiation 
of the infinitesimal action of its Lie algebra $\hbar(\al\otimes S\hl)^\hl[[\hbar]]$: 
\begin{equation}\label{eq:infinit}
q\cdot\alpha={\rm d}q+[q,\alpha]\qquad\big(q\in\hbar(\al\otimes S\hl)^\hl[[\hbar]]\big)
\end{equation}

Then going along the lines of subsection \ref{sec:p2} one can prove the 
following 
\begin{prop}
There exists a $L_\infty$-quasi-isomorphism 
$$\mathcal Q_1:(C^c((\wedge^*\ml)^\hl[2]),[,]_{\ml})\to(\ggun,{\rm d}+[,])$$
such that $\mathcal Q_1^1$ is the natural inclusion and $\mathcal Q_1^n$ takes values in $(\gun)_{\leq n-1}$. 
\end{prop}
Consequently any Maurer-Cartan element in $(\gun[1])[[\hbar]]$ is equivalent to a one of the form 
$\hbar\rho_\hbar(\hbar\lambda)$, where $\rho_\hbar\in(\wedge^2\al\hat\otimes\hat S\hl)^\hl[[\hbar]]$ 
satisfies CDYBE. In other words $\rho_\hbar$ is $\hbar$-dependant formal dynamical $r$-matrix. 
On such a $\rho_\hbar$ the infinitesimal action (\ref{eq:infinit}) becomes 
\begin{equation}\label{gaugermat}
q\cdot\rho_\hbar=-\sum_ih_i\wedge\frac{\partial q}{\partial \lambda^i}+[q,\rho_\hbar]
\qquad(q\in\al\hat\otimes\hat S\hl)^\hl[[\hbar]]
\end{equation}
This action integrates in an affine action of some group $\widetilde{G_1}$ of $\hl$-equivariant formal 
maps with values in the Lie group $G$ of $\al$. And then we have a bijection 
\begin{equation}\label{classifrmat}
\frac{\{\pi\in\hbar(\wedge^2\ml)^\hl[[\hbar]]~\textrm{s.t.}~[\pi,\pi]_\ml=0\}}{G_0}
\longleftrightarrow
\frac{\{\textrm{form. dynam. $r$-matrices}/\R[[\hbar]]\}}{\widetilde{G_1}}
\end{equation}
\begin{rmq}\emph{
This bijection has to be compared with Proposition 2.13 in \cite{X1} and section 3 of \cite{ES}
}\end{rmq}

Finally, combining (\ref{classifrmat}), (\ref{eq:classiftwist}) and (\ref{eq:algform}) we obtain the 
following generalization of Theorem 6.11 in \cite{X1} to the case of a nonabelian base: 
\begin{thm}\label{thm:class}
Let $\pi\in(\wedge^2\ml)^\hl$ such that $[\pi,\pi]_\ml=0$. Then there are bijective correspondances 
between 
\begin{enumerate}
\item the set of $\hbar$-dependant and $G$-invariant Poisson structures $\pi_\hbar=\hbar\pi~\textrm{mod}~\hbar^2$ 
on $G/H$, modulo the action of $G_0$, 
\item the set of $\hbar$-dependant formal dynamical $r$-matrices $\rho_\hbar(\lambda)$ such that 
$\rho_\hbar(0)=\pi~\textrm{mod}~\hbar$ in $\wedge^2(\al/\hl)[[\hbar]]$, 
modulo the action (\ref{gaugermat}) of $\widetilde{G_1}$, 
\item the set of formal dynamical twists $J(\lambda)$ satisfying 
$\textrm{Alt}\frac{J(0)-1}{\hbar}=\pi~\textrm{mod}~\hbar$ in $\wedge^2(\al/\hl)[[\hbar]]$, 
modulo gauge equivalence (\ref{eq:gaugequiv}). 
\end{enumerate}
\end{thm}

\section{Another case when the twist quantization exists}\label{sec:ex}

In this section we assume that $\hl$ is abelian and admits a Lie subalgebra $\ml$ as complement. 

Note that since $\hl$ is abelian and $\ml$ a Lie subalgebra, the projection $p:\al\to\al$ on $\ml$ along $\hl$ extends to 
a morphsim of graded Lie algebras $\wedge^{\cdot}p:(\wedge\al)^\hl\to(\wedge\al)^\hl$ at the level of $\hl$-invariants. 
And thus $\wedge^{\cdot}p\otimes\varepsilon:(\gun[1],{\rm d},[,])\to((\wedge\al)^\hl[1],0,[,])$ is a morphism of 
{\tt dgla}'s. 
Then the natural inclusion ${\rm id}\otimes1:(T^*U\al)^\hl\to\gde$ obviously allows one to consider $(T^*U\al)^\hl[1]$ as 
a sub-{\tt dgla} of $\gde[1]$. Finally recall from \cite[section 3.3]{Ca} that there exists a $L_\infty$-quasi-isomorphism 
$\mathcal F:C^c((\wedge^*\al)^\hl[2])\to C^c((T^*U\al)^\hl[2])$ with $\mathcal F^1={\rm alt}$. By composing these maps one 
obtains a $L_\infty$-morphism $$\widetilde{\mathcal F}:(\Gun,{\rm d}+[,])\to(\Gde,b+[,]_G)$$ with values in 
$(\Gde)_{\leq0}$ and first structure map 
$\widetilde{\mathcal F}^1=({\rm alt}\otimes1)\circ(\wedge^{\cdot}p\otimes\varepsilon)$. 
\begin{thm}
There exists a $L_\infty$-quasi-iomorphism 
$$\Psi:(\Gun,{\rm d}+[,])\to(\Gde,b+[,]_G)$$
with properties {\rm (F1)} and {\rm (F2)} of Theorem \ref{thm:main}. 
\end{thm}
\begin{proof}
First observe that since $\hl$ is abelian then $\gun\cong((\wedge\al)^\hl\cap\wedge\ml)\otimes\wedge\hl\otimes S\hl$ as a 
vector space. Thus if $\delta_K$ is as in the proof of proposition \ref{reso1} then $\delta:={\rm id}\otimes\delta_K$ is a 
homotopy operator: $\delta{\rm d}+{\rm d}\delta={\rm id}-\wedge^{\cdot}p\otimes\varepsilon$ and $\delta\circ\delta=0$. 

Now we proceed like in subsection \ref{sec:end}: use lemma \ref{lem} to construct a $L_\infty$-morphism 
$\Psi$ with first structure map $\Psi^1=\widetilde{\mathcal F}^1+b\circ V+V\circ{\rm d}$, where 
$V:=({\rm alt}\otimes{\rm sym})\circ\delta:\gun\to\gde[-1]$. 

It remains to prove that $\Psi$ is a quasi-isomorphism. 
It follows from the first observation in this proof that $H^*(\gun,{\rm d})=(\wedge\al)^\hl\cap\wedge\ml$, which also 
equals $H^*(\gde,b)$ due to the spectral sequence argument. Consequently $\widetilde{\mathcal F}^1$ is a quasi-isomorphism 
of complexes, and so is $\Psi^1$. 
\end{proof}
Finally using the same argumentation as in the proof of theorem \ref{thm:rmat} (subsection \ref{sec:proof}) one obtains the 
\begin{thm}
If $\hl$ is an abelian subalgebra of $\al$ with a Lie subalgebra as a complement, then any formal classical dynamical 
$r$-matrix for $(\al,\hl,0)$ admits a dynamical twist quantization (associated to the trivial associator). 
\end{thm}
\begin{ex}\emph{
In particular, this allows us to quantize dynamical $r$-matrices arizing from semi-direct products 
$\al=\ml\ltimes\C^n$ like in \cite[example 3.7]{EN}. 
}\end{ex}

\section*{Concluding remarks}

Let us first observe that if $\hl$ is abelian then $(\wedge^*\al)^\hl\cap\wedge^*\ml[1]$ 
({\rm resp.} $(T^*U\al)^\hl\cap T^*{\rm sym}(S\ml)[1]$) inherits a {\tt dgla} structure from the one of $\gun[1]$ 
({\rm resp.} $\gde[1]$) and $H^*(\gun,{\rm d})=(\wedge^*\al)^\hl\cap\wedge^*\ml=H^*(\gde,b)$, for any complement $\ml$ 
of $\hl$. 
Thus I conjecture that there exists a $L_\infty$-quasi-isomorphism between 
$(\wedge^*\al)^\hl\cap\wedge^*\ml[1]$ and $(T^*U\al)^\hl\cap T^*{\rm sym}(S\ml)[1]$ which generalizes together $\phi$ of 
subsection \ref{sec:end} and $\mathcal F$ of section \ref{sec:ex}. In particular this would imply conjecture \ref{thm:cjt} 
in the abelian (and non-modified) case. \\

Let us then mention that one can consider a non-triangular (i.e., non-antisymmetric) version of non-modified classical 
dynamical $r$-matrices. Namely, $\hl$-equivariant maps $r\in{\rm Reg}(\hl^*,\al\otimes\al)$ such that 
$\textrm{CYB}(r)-\textrm{Alt}(dr)=0$. According to \cite{X}, a quantization of such a $r$ is a $\hl$-equivariant 
map $R=1+\hbar r+O(\hbar^2)\in{\rm Reg}(\hl^*,U\al^{\otimes2})[[\hbar]]$ that satisfies the 
\emph{quantum dynamical Yang-Baxter equation} (QDYBE) 
\begin{equation}
R^{1,2}(\lambda)*R^{1,3}(\lambda+\hbar h^2)*R^{2,3}(\lambda)
=R^{2,3}(\lambda+\hbar h^1)*R^{1,3}(\lambda)*R^{1,2}(\lambda+\hbar h^3)
\end{equation}
\begin{qtn}
Does such a quantization always exist?
\end{qtn}
The most famous example of non-triangular dynamical $r$-matrices was found in \cite{AM} by Alekseev and Meinrenken, then 
extended successively to a more general context in \cite{EV,ES,EE1}, and quantized in \cite{EE1}. \\

Following \cite{EE1}, remark that for any non-triangular dynamical $r$-matrix $r$ such that 
$r+r^{{\rm op}}=t\in(S^2\al)^\al$ (quasi-triangular case) one can define $\rho:=r-t/2$ and $Z:=\frac14[t^{1,2},t^{2,3}]$. 
Then $\rho$ is a modified dynamical $r$-matrix for $(\al,\hl,Z)$; morever the assignment $r\longmapsto\rho$ is a 
bijective correspondence between  quasi-triangular dynamical $r$-matrices for $(\al,\hl,t)$ and modified dynamical 
$r$-matrices for $(\al,\hl,Z)$. Now observe that if $J(\lambda)$ is a dynamical twist quantizing $\rho$, then 
$R(\lambda)=J^{{\rm op}}(\lambda)^{-1}*e^{\hbar t/2}*J(\lambda)$ is a quantum dynamical $R$-matrix quantizing $r$. 

In this paper we have constructed such a dynamical twist in the triangular case $t=0$. One can ask 
\begin{qtn}
Does such a dynamical twist exist for any quasi-triangular dynamical $r$-matrix? At least in the reductive and abelian 
cases? 
\end{qtn}
This question seems to be more reasonable than the previous one. \\
More generally one can ask if conjecture \ref{thm:cjt} (and its smooth and meromorphic versions) is true in general. A 
positive answer was given in \cite{EE1} when $\hl=\al$; but unfortunately it is not known in general, even for the 
non-dynamical case $\hl=\{0\}$ (which is the last problem of Drinfeld \cite{D}: quantization of coboundary Lie 
bialgebras). \\

Finally let me mention that if $r(\lambda)$ is a triangular dynamical $r$-matrix for $(\al,\hl)$, then the bivector field 
$$\pi:=\overrightarrow{r(\lambda)}+\sum_i\frac\partial{\partial\lambda^i}\wedge\overrightarrow{h_i}+\pi_{\hl^*}$$
is a $G\times H$-biinvariant Poisson structure on $G\times\hl^*$ and the projection $p:G\times\hl^*\rightarrow\hl^*$ is a 
momentum map. Moreover, according to \cite{X} any dynamical twist quantization $J(\lambda)$ of $r(\lambda)$ allows us to 
define a $G\times H$-biinvariant star-product $*$ quantizing $\pi$ on $G\times\hl^*$ as follows: 
\begin{equation*}
\begin{array}{lll}
f*g=f*_{PBW}g & \textrm{if} \quad & f,\,g\in C^\infty(\hl^*) \\
f*g=fg & \textrm{if} \quad & f\in C^\infty(G),\,g\in C^\infty(\hl^*) \\
f*g=\exp\big(\hbar\sum_i\frac\partial{\partial\lambda^i}\otimes\overrightarrow{h_i}\big)\cdot(f\otimes g) & \textrm{if} 
\quad & f\in C^\infty(\hl^*),\,g\in C^\infty(G) \\
f*g=\overrightarrow{J(\lambda)}(f\otimes g) & \textrm{if} \quad & f,\,g\in C^\infty(G)
\end{array}
\end{equation*}
This way the map $p^*:({\rm Fct}(\hl^*)[[\hbar]],*_{PBW})\to({\rm Fct}(G\times\hl^*)[[\hbar]],*)$ becomes a quantum 
momentum map in the sens of \cite{X2}. \\
So there may be a way to see momentum maps and their quantum analogues as Maurer-Cartan elements in {\tt dgla}'s. \\

\begin{appendix}

\section{Homotopy Lie algebras}

See \cite{HS} for a detailed discussion of the theory. \\

Recall that a \emph{$L_\infty$-algebra} structure on a graded vector space $\al$ is a degree $1$ coderivation $Q$ on the 
cofree cocommutative coalgebra $C^c(\al[1])$ such that $Q\circ Q=0$. By cofreeness, such a coderivation $Q$ is uniquely 
determined by structure maps $Q^n:\Lambda^n\al\to\al[2-n]$ which satisfy an infinite collection of equations. In 
particular $(\al,Q^1)$ is a cochain complex. 
\begin{ex}\emph{
Any {\tt dgla} $(\al,{\rm d},[,])$ is canonically a $L_\infty$-algebra. Namely, $Q$ is given by structure maps 
$Q^1={\rm d}$, $Q^2=[,]$ and $Q^n=0$ for $n>2$. 
}\end{ex}

A \emph{$L_\infty$-morphism} between two $L_\infty$-algebras $(\gun,Q_1)$ and $(\gde,Q_2)$ is a degree $0$ morphism of 
coalgebras $F:C^c(\gun[1])\to C^c(\gde[1])$ such that $F\circ Q_1=Q_2\circ F$. Again by cofreeness, such a 
morphism is uniquely determined by structure maps $F^n:\Lambda^n\gun\to\gde[1-n]$ which satisfy an infinite 
collection of equations. In particular $F^1:\gun\to\gde$ is a morphism of complexes; when it induces an isomorphism in 
cohomology we say that $F$ is a \emph{$L_\infty$-quasi-isomorphism}. 
\begin{ex}\emph{
Any morphism of {\tt dgla}'s is a $L_\infty$-morphism with all structure maps equal to zero except the first one. 
}\end{ex}

In this paper we use many times the following 
\begin{lem}[\cite{Do}]\label{lem}
Let $F:C^c(\gun[1])\to C^c(\gde[1])$ be a $L_\infty$-morphism. For any linear map $V:\gun\to\gde[-1]$ 
there exists a $L_\infty$-morphism $\Psi:C^c(\gun[1])\to C^c(\gde[1])$ with first structure map 
$\Psi^1=F^1+Q_2^1\circ V+V\circ Q_1^1$. Moreover, if $F$ is a $L_\infty$-quasi-isomorphism then $\Psi$ is also. 
\end{lem}
\begin{proof}
First remark that $V$ extends uniquely to a linear map $C^c(\gun[1])\to C^c(\gde[1])$ of degree $-1$ such that 
$$\Delta_2\circ V=\big(F\otimes V+V\otimes F+
\frac12V\otimes(Q_2\circ V+V\circ Q_1)+\frac12(Q_2\circ V+V\circ Q_1)\otimes V\big)\circ\Delta_1$$
where $\Delta_1$ and $\Delta_2$ denote comultiplications in $C^c(\gun[1])$ and $C^c(\gde[1])$, respectively. \\
Then define $\Psi:=F+Q_2\circ V+V\circ Q_1$. 
\end{proof}
\begin{rmq}\label{remark}\emph{
Assume that in the previous lemma $\gun$ and $\gde$ are filtrated, $F$ is such that $F^n$ takes values in 
$(\gde)_{\leq n-1}$, and $V\big((\gun)_{\leq k}\big)\subset(\gde)_{\leq k+1}$. Then one can obviously check that 
for any $X\in(\Lambda^n\gun)_{\leq k}$, $F^n(X)\in(\gde)_{\leq n+k-1}$. }
\end{rmq}

\end{appendix}

\mbox{}

\noindent\footnotesize{\textsc{IRMA, 7 rue René Descartes, F-67084 Strasbourg, France} \\
\emph{E-mail address}: {\bf calaque@math.u-strasbg.fr}


\thebibliography{42}

\bibitem[AM]{AM}
A.~Alekseev, E.~Meinrenken, The non-commutative Weil algebra, Invent. Math. {\bf 139} (2000), 135-172. 

\bibitem[Ba]{B}
J.~H.~Baues, The double bar and cobar construction, Compos. Math. {\bf 43} (1981), 331-341. 

\bibitem[Ca]{Ca}
D.~Calaque, Formality for Lie algebroids, to appear in Comm. Math. Phys. 

\bibitem[Do]{Do}
V.~Dolgushev, Covariant and equivariant formality theorem, Adv. Math. {\bf 191} (2005), 147-177. 

\bibitem[Dr1]{D}
V.~Drinfeld, On some unsolved problems in quantum group theory, Lecture Notes Math. {\bf 1510} (1992), 1-8. 

\bibitem[Dr2]{D2}
V.~Drinfeld, Quasi-Hopf algebras, Leningrad Math. J. {\bf 1} (1990), 1419-1457. 

\bibitem[EE1]{EE1}
B.~Enriquez, P.~Etingof, Quantization of Alekseev-Meinrenken dynamical $r$-matrices, Trans. Am. Math. Soc. (ser. 2) 
{\bf 210} (2003), 81-98. 

\bibitem[EE2]{EE2}
B.~Enriquez, P.~Etingof, Quantization of classical dynamical $r$-matrices with nonabelian base, to appear in Comm. Math. 
Phys.

\bibitem[EN]{EN}
P.~Etingof, D.~Nikshych, Vertex-IRF transformations and quantization of dynamical $r$-matrices, Math. Res. Lett. {\bf 8} 
(2001), 331-345. 

\bibitem[ES]{ES}
P.~Etingof, O.~Schiffmann, On the moduli space of classical dynamical $r$-matrices, Math. Res. Lett. {\bf 8} (2001), 
157-170. 

\bibitem[EV]{EV}
P.~Etingof, A.~Varchenko, Geometry and classification of solutions of the classical dynamical Yang-Baxter equation, 
Comm. Math. Phys. {\bf 192} (1998), 77-120. 

\bibitem[Fe]{F}
G.~Felder, Conformal field theory and integrable systems associated to elliptic curves, Proceedings of the International 
Congress of Mathematicians, Vol. 1, 2 (Zurich, 1994), 1247-1255, Birkhäuser, Basel, 1995. 

\bibitem[Ge]{G}
E.~Getzler, Cartan homotopy formula and the Gauss-Manin connection in cyclic homology, Israel Math. Conf. Proc. {\bf 102} 
(1993), 256-283. 

\bibitem[He]{H}
S.~Helgason, \emph{Groups and geometric analysis}, Pure Appl. Math., Vol. 113, Orlando, 1984. 

\bibitem[HS]{HS}
V.~Hinich, V.~Schechtman, Homotopy Lie algebras, I.M.~Gelfand Seminar, Adv. Sov. Math. {\bf 16} (1993), no. 2, 1-28. 

\bibitem[Kh]{Kh}
M.~Khalkhali, Operations on cyclic homology, the $X$ complex, and a conjecture of Deligne, Comm. Math. Phys. {\bf 202} 
(1999), no. 2, 309-323. 

\bibitem[Ko]{K}
M.~Kontsevich, Deformation quantization of Poisson manifolds, Lett. Math. Phys. {\bf 66} (2003), no. 3, 157-216. 

\bibitem[Mo]{Mo}
T.~Mochizuki, An application of formality theorem to a quantization of dynamical $r$-matrices, unpublished preprint. 

\bibitem[No]{N}
K.~Nomizu, Invariant affine connections on homogeneous spaces, Amer. J. Math. {\bf 76} (1954), 33-65. 

\bibitem[Ta]{T}
D.~Tamarkin, Another proof of M.~Kontsevich formality theorem for $\mathbb R^n$, preprint math.QA/9803025 (1998). 

\bibitem[X1]{X2}
P.~Xu, Fedosov $*$-products and quantum momentum maps, Comm. Math. Phys. {\bf 197} (1998), 167-197. 

\bibitem[X2]{X1}
P.~Xu, Triangular dynamical $r$-matrices and quantization, Adv. Math. {\bf 166} (2002), no. 1, 1-49. 

\bibitem[X3]{X}
P.~Xu, Quantum dynamical Yang-Baxter equation over a nonabelian base, Comm. Math. Phys. {\bf 226} (2002), no. 3, 475-495.
\end{document}